# NONANTICIPATING ESTIMATION APPLIED TO SEQUENTIAL ANALYSIS AND CHANGEPOINT DETECTION

By Gary Lorden[1] and Moshe Pollak[2]

*California Institute of Technology and Hebrew University of Jerusalem*

Suppose a process yields independent observations whose distributions belong to a family parameterized by $\theta \in \Theta$. When the process is in control, the observations are i.i.d. with a known parameter value $\theta_0$. When the process is out of control, the parameter changes. We apply an idea of Robbins and Siegmund [*Proc. Sixth Berkeley Symp. Math. Statist. Probab.* **4** (1972) 37–41] to construct a class of sequential tests and detection schemes whereby the unknown post-change parameters are estimated. This approach is especially useful in situations where the parametric space is intricate and mixture-type rules are operationally or conceptually difficult to formulate. We exemplify our approach by applying it to the problem of detecting a change in the shape parameter of a Gamma distribution, in both a univariate and a multivariate setting.

**1. Introduction.** In all but the simplest cases, the problem of detecting a change involves at least one unknown post-change parameter. In the well-known Shiryayev–Roberts detection scheme [11, 12], a change from parameter value $\theta = \theta_0$ (possibly multidimensional) to $\theta = \theta_1$, say, in the distribution of a sequence of i.i.d. observations $X_1, X_2, \ldots$ is detected by a stopping rule

$$N_A = \min\{n | R_n \geq A\},$$

Received February 2000; revised March 2004.
[1]Supported in part by a Lady Davis Fellowship at the Hebrew University of Jerusalem.
[2]Supported by a grant from the Israel Science Foundation. Also supported by the Marcy Bogen Chair of Statistics at the Hebrew University of Jerusalem.

where
$$R_n = \sum_{k=1,\ldots,n} \Lambda_{n,k}$$
and
$$\Lambda_{n,k} = \prod_{i=k,\ldots,n} f_{\theta_1}(X_i)/f_{\theta_0}(X_i).$$

When the post-change $\theta_1$ is not known and it is desired to respond quickly to a broad range of possible values, the Shiryayev–Roberts (SR) rule is in principle easy to modify: just introduce a mixing measure $\lambda(\theta)$ and define

$$\Lambda_{n,k} = \int_\Theta \left( \prod_{i=k,\ldots,n} (f_\theta(X_i)/f_{\theta_0}(X_i)) \right) \lambda(\theta)\, d\theta.$$

As is well known [5], this approach preserves the martingale property of the sequence $\{R_n - n\}$ under the "no change" probability measure, $P_\infty$, so that

$$E_\infty N_A = E_\infty R_{N_A} \geq A,$$

a useful lower bound on the average run length (ARL) to false alarm. Moreover, it is typically true that [5]

$$E_\infty N_A / A \to 1/\gamma \qquad \text{as } A \to \infty,$$

where $\gamma$ can be either evaluated by renewal-theoretic methods or simulated, which suggests using the approximation

$$E_\infty N_A \approx A/\gamma.$$

In practice, however, it is usually difficult to carry out the computation of $\Lambda_{n,k}$'s unless the mixing measure $\lambda$ can be chosen as a natural conjugate prior. Moreover, in many cases, particularly when $\theta$ is multidimensional, it is conceptually difficult to make natural choices of $\lambda$. The present paper suggests an alternative approach, based upon defining

$$\Lambda_{n,k} = \prod_{i=k,\ldots,n} f_{\theta_{i,k}}(X_i)/f_{\theta_0}(X_i),$$

where $\theta_{i,k}$ is an estimator of $\theta$ based upon $X_k,\ldots,X_{i-1}$. The same idea appears in [9] in the context of sequential hypothesis testing and is applied to the changepoint problem in [1, 3]. By requiring $\theta_{n,k}$ not to depend on $X_n$, one preserves the martingale property of $\{R_n - n\}$ and similarly the upper bound on significance levels that Robbins and Siegmund rely upon. The potential advantage of this approach in complicated settings is that simple estimators based on the method of moments or maximum likelihood are usually much easier to choose than appropriate mixtures, as well as substantially faster to compute.

Moreover, in many cases:

1. The asymptotic "overshoot correction," $\gamma$, is finite and can be evaluated readily by simulating Robbins–Siegmund-type hypothesis tests rather than the changepoint detection rules themselves.
2. The proposed Shiryayev–Roberts–Robbins–Siegmund (SRRS) detection rules have reasonably good efficiency, that is, short post-change delays to detection, nearly as good as mixture rules.

Since the overshoot constant $\gamma$ is most easily evaluated in the context of hypothesis testing, Section 2 and part of Section 3 are devoted to problems of testing. Proofs of the asymptotic analysis of the operating characteristics of the SRRS rules involve formidable calculations. Therefore, our approach in the present paper is to illustrate the arguments in special cases, the testing of hypotheses about the mean of a normal distribution (Section 2) and hypothesis testing and changepoint detection for the shape parameter, $\theta$, of a Gamma distribution (Sections 3–5). We believe that these special cases provide a good introduction to the type of argument that will work in many other contexts.

Section 2 illustrates the pattern of the asymptotic behavior of the estimator sequence and the determination of $\gamma$. It turns out that there is a natural correspondence between choices of an estimator sequence and a choice of mixture $\lambda$, suggesting that, at least asymptotically, the two methods have a natural relationship. In Sections 3 and 5 we give asymptotic results for the Gamma shape example, showing in particular that the asymptotic efficiency of the estimator sequence used in the SRRS scheme determines the coefficient of the second-order term in the asymptotic expansion for the expected delay to detection. In particular, an asymptotically efficient sequence of estimators yields a second-order asymptotically optimal detection scheme. (For comparison, Dragalin's [1] scheme does not achieve this.) Section 4 gives results of Monte Carlo simulations of the performance of the SRRS scheme using the method-of-moments and maximum likelihood estimators of the Gamma shape parameter, as well as comparisons with other changepoint detection rules. Section 5 illustrates the application of the SRRS approach to multiparameter problems.

**2. A first example: tests for the value of a normal mean.** Let $X_1, X_2, \ldots$ be a sequence of independent $N(\mu, 1)$-distributed random variables, and suppose one is interested in a power one $\alpha$-level test of $H_0 : \mu = 0$ versus $H_1 : \mu \neq 0$. Robbins and Siegmund [9] introduced the following sequential test: Let $\mu_n$ be an $F_{n-1} = F(X_1, \ldots, X_{n-1})$-measurable estimate of $\mu$ (where $F_0$ is the trivial $\sigma$-field), define $\Lambda_n = \exp\{\sum_{i=1,\ldots,n}(\mu_i X_i - (\mu_i)^2/2)\}$, $\tau_b = \min\{n | n \geq 1, \Lambda_n \geq \exp\{b\}\}$ ($\tau_b = \infty$ if no such $n$ exists); reject $H_0$ if and only if $\tau_b < \infty$. By using the martingale property of $\Lambda_n$ under $H_0$, Robbins and Siegmund showed that $\alpha = P_{H_0}(\tau_b < \infty) \leq \exp\{-b\}$. In this section, we will give an approximation for $\alpha$ for a special case of the sequence $\{\mu_n\}$ when $b$ is large.

Let $\mu_1 = s/t$ ($= 0$ if $s = t = 0$) and $\mu_{n+1} = (n\bar{X}_n + s)/(n+t)$, where $\bar{X}_n = \sum_{i=1,\ldots,n} X_i/n$ and $-\infty < s < \infty$, $0 < t < \infty$ or $s = t = 0$ are constants. (This would be a natural estimate of $\mu$ after $n$ observations on test when prior to testing there is a learning sample of $t$ observations whose sum is $s$; it is also a way of incorporating a prior distribution into the testing scheme.) In other words, after every observation we update our estimate of $\mu$, and $\exp\{\mu_n X_n - (\mu_n)^2/2\}$ is the estimated likelihood ratio for the $n$th observation.

Let $G_{s,t}$ be the $N(s/t, \sum_{i=1,\ldots,\infty} 1/(i+t)^2)$ c.d.f. (where $s/t = 0$ if $s = t = 0$). Let $\mu_1 = s/t$, $v^2(t) = \sum_{i=1,\ldots,\infty} 1/(i+t)^2$ and

$$\nu(\mu) = 2\mu^{-2} \exp\left\{-2 \sum_{n=1,\ldots,\infty} n^{-1} \Phi(-\tfrac{1}{2}|\mu|\sqrt{n})\right\}. \tag{1}$$

THEOREM 1. *As $b \to \infty$, when $\mu_1 = s/t$ ($= 0$ if $s = t = 0$) and $\mu_{n+1} = (n\bar{X}_n + s)/(n+t)$, the Robbins–Siegmund power one test of $H_0 : \mu = 0$ versus $H_1 : \mu \neq 0$ has significance level*

$$\alpha = P_{H_0}(\tau_b < \infty) = (1 + o(1))\gamma \exp\{-b\},$$

$$\text{where } \gamma = \int_{-\infty}^{\infty} \nu(y)\, dG_{s,t}(y). \tag{2}$$

REMARK. A theorem analogous to Theorem 1 can be formulated for an arbitrary sequence $\{\mu_n\}$. Its practical value is usually as a statement of existence of a limit, which one can evaluate by simulation. The analog to $G_{s,t}$ is generally very hard to compute.

PROOF OF THEOREM 1. Let $Q$ be the measure on $\{X_1, X_2, \ldots\}$ under which the distribution of $X_n$ conditional on $X_1, \ldots, X_{n-1}$ is $N(\mu_n, 1)$; $n = 1, 2, \ldots$. [By abuse of notation, we will let $Q(X_1, \ldots, X_n)$ denote the distribution of $X_1, \ldots, X_n$ under $Q$.] Let $P_Q$ and $E_Q$ denote probability and expectation, respectively, under the measure $Q$. The proof requires two lemmas.

LEMMA 1. *Under the measure $Q$, the sequence $\{\mu_n\}$ converges a.s. to a random variable whose distribution is $G_{s,t}$.*

PROOF. Write $X_n = \mu_n + Z_n$ where $Z_i \sim N(0, 1)$ and are independent. Thus, for $n \geq 2$

$$\mu_{n+1} = \left(\sum_{i=1,\ldots,n} X_i + s\right) \Big/ (n+t)$$

$$= ((n-1+t)\mu_n + \mu_n + Z_n)/(n+t) = \mu_n + Z_n/(n+t).$$

Therefore
$$\mu_n = \mu_1 + \sum_{i=1,\ldots,n-1} Z_i/(i+t).$$

Hence $\mu_n$ converges a.s. as $n \to \infty$ to $\mu_1 + \sum_{i=1,\ldots,\infty} Z_i/(i+t)$, whose distribution is $G_{s,t}$. □

LEMMA 2. $P_Q(\tau_b < \infty) = 1$.

PROOF. By virtue of Lemma 1, $\sum_{i=1,\ldots,n}(\mu_i X_i - (\mu_i)^2/2) \to \infty$ a.s. $(Q)$ as $n \to \infty$, from which Lemma 2 follows. □

PROOF OF THEOREM 1—CONTINUED. Let $\Lambda_n = \Lambda_n(X_1,\ldots,X_n) = \exp\{\sum_{i=1,\ldots,n}(\mu_i X_i - (\mu_i)^2/2)\}$. Since $P_Q(\tau_b < \infty) = 1$,

$$\alpha = P_{H_0}(\tau_b < \infty) = \sum_{n=1}^{\infty} \int_{\tau_b=n} \cdots \int f_{H_0}(x_1,\ldots,x_n)\,dx_1\cdots dx_n$$

$$(3) \quad = \exp\{-b\} \sum_{n=1}^{\infty} \int_{\tau_b=n} \cdots \int \exp\{-(\log \Lambda_n - b)\}\,dP_Q(X_1,\ldots,X_n)$$

$$= \exp\{-b\} E_Q \exp\{-(\log \Lambda_{\tau_b} - b)\}.$$

Thus what is left to be done is a renewal-theoretic analysis of the expectation in (3).

Let $0 < \varepsilon < 1$. By virtue of Lemma 1, there exist $0 < d_\varepsilon < a_\varepsilon < \infty$ such that $P_Q(d_\varepsilon < |\mu_n| < a_\varepsilon \text{ for all } n \geq 2) \geq 1 - \varepsilon$. Note that if $Y \sim N(\mu,1)$, then $\sup_{a>0} P(Y - a > y|Y > a) \to 0$ as $y \to \infty$ uniformly for all $\mu$ in a compact set. Therefore, there exists $0 < c_\varepsilon < \infty$ such that if $b > c > c_\varepsilon$, then $P_Q(A_{b,c}) \geq 1 - 2\varepsilon$, where $A_{b,c} = \{\log \Lambda_{\tau_{b-c}} - (b-c) \leq c/2\}$.

Choose $c > c_\varepsilon$ and write $w = \tau_{b-c}$. Note that $X_{w+j} = \mu_w + \sum_{i=0,\ldots,j-1} Z_{w+i}/(w+i+t) + Z_{w+j}$. When $j$ remains fixed and $b \to \infty$, then $\sum_{i=1,\ldots,j-1} Z_{w+i}/(w+i) \to 0$ in probability. Leaving $c$ fixed, when $\mu_w \notin (-d_\varepsilon, d_\varepsilon)$, $\tau_b - w = \tau_b - \tau_{b-c}$ is stochastically bounded in $Q$-probability as $b \to \infty$.

For $\mu \neq 0$, let $H_{\mu,b} = \min\{n \geq 1, \sum_{i=1,\ldots,n}(\mu(Z_i + \mu) - \mu^2/2) \geq b\}$. Note that normally distributed variables are strongly nonlattice, so that the convergence in distribution of $\sum_{i=1,\ldots,H_{\mu,b}}(\mu(Z_i + \mu) - \mu^2/2) - b$ to its renewal-theoretic limit is uniform on compact sets that do not contain zero (see [16]). Therefore, for large enough $c$, for all $d_\varepsilon < |\mu| < a_\varepsilon$, on $\sum_{i=1,\ldots,H_{\mu,b-c}}(\mu(Z_i + \mu) - \mu^2/2) < b - c/2\}$,

$$(4) \quad \left| E\left(\exp\left\{-\left(\sum_{i=1,\ldots,H_{\mu,b}}(\mu(Z_i+\mu) - \mu^2/2) - b\right)\right\}\bigg| F_{H_{\mu,b-c}}\right) - \nu(\mu) \right| < \varepsilon$$

($\nu(\mu)$ of (1) is the renewal-theoretic limit of the expectation; cf. [13]).

Also, for fixed $k$

$$\max_{1\leq j\leq k}\left|\sum_{i=w+1,\ldots,w+j}(\mu_i(X_i+\mu_i)-\mu_i^2/2)\right.$$

$$\left.-\sum_{i=w+1,\ldots,w+j}(\mu_w(X_i+\mu_w)-\mu_w^2/2)\right|\xrightarrow{P}0 \quad \text{as } b\to\infty.$$

Therefore, and because of the stochastic boundedness of $\tau_b - w$, for all large enough $b$, on $A_{b,c}\cap\{d_\varepsilon<|\mu_w|<a_\varepsilon\}$

$$\left|E_Q(\exp\{-(\log\Lambda_{\tau_b}-b)\}|F_w;\mu_w=y;\log\Lambda_w=z)\right.$$

(5)
$$-E\left(\exp\left\{-\left(\sum_{i=1,\ldots,H_{y,b}}(y(Z_i+y)-y^2/2)-b\right)\right\}\right.$$

$$\left.\left.\left|\sum_{i=1,\ldots,w}(y(Z_i+y)-y^2/2)=z\right)\right|<\varepsilon.$$

Combining (4) and (5), one obtains that $c$ can be fixed so that for all large enough $b$, on $A_{b,c}\cap\{d_\varepsilon<|\mu_w|<a_\varepsilon\}$

$$|E_Q(\exp\{-(\log\Lambda_{\tau_b}-b)\}|F_w)-\nu(\mu_w)|<2\varepsilon.$$

Since the distribution of $\mu_w$ converges as $b\to\infty$ to $G_{s,t}$, it follows that one can fix $c$ so that there exists $b_{\varepsilon,c}$ such that for all $b>b_{\varepsilon,c}$

$$\left|E_Q\left(\exp\{-(\log\Lambda_{\tau_b}-b)\}\right)-\int_{-\infty}^{\infty}\nu(y)\,dG_{s,t}(y)\right|<6\varepsilon.$$

Letting $\varepsilon\to 0$ concludes the proof of Theorem 1. $\square$

As the preceding analysis suggests, every stopping rule of the Robbins–Siegmund type has a mixture analog. For example, the mixture-type analog of the rule described in Theorem 1 is

$$T_b=\min\left\{n\left|\int\exp\left\{y\sum_{i=1,\ldots,n}X_i-ny^2/2\right\}dG_{s,t}(y)\geq\exp\{b\}\right.\right\}$$

$$=\min\left\{n\left|\exp\left\{\left(s/t+v^2(t)\sum_{i=1,\ldots,n}X_i\right)^2\Big/[2v^2(t)(nv^2(t)+1)]\right\}\right.\right.$$

$$\left.\left.\times(nv^2(t)+1)^{-1/2}\geq\exp\{b\}\right\}\right.$$

and [2] the asymptotic expression for its level of significance $P_{H_0}(T_b < \infty)$ is the same as (2). For a given level of significance $\alpha$, both rules have (approximately) the same stopping threshold $\exp\{b\}$.

Robbins and Siegmund [10] noted that in the continuous-time (Brownian motion) case the two rules are identical. In the discrete-time case a comparison between them is of interest. Following the methods of Pollak and Siegmund [6], Robbins and Siegmund [10], Lai and Siegmund [2] and Woodroofe [17], one can show that, for any fixed $\mu \neq 0$, the difference between the expected sample sizes of the two stopping rules is $O(1)$ as $\alpha \to 0$. Specifically, letting $r_t = \lim_{n \to \infty}(\sum_{i=1,\ldots,n} 1/(i+t) - \log n)$, it can be shown that

$$E_\mu \tau_b - E_\mu T_b$$
$$= \mu^{-2}\Bigg\{\Bigg[r_t - t \sum_{j=1,\ldots,\infty} 1/(j+t)^2 - 2\log\Bigg(\sum_{j=1,\ldots,\infty} 1/(j+t)^2\Bigg) + 1\Bigg]$$
$$\qquad + (\mu - s/t)^2\Bigg[1 + t^2 \sum_{j=1,\ldots,\infty} 1/(j+t)^2 - 1\Big/\sum_{j=1,\ldots,\infty} 1/(j+t)^2\Bigg]\Bigg\}$$
$$\qquad + o(1)$$
$$\stackrel{\text{def}}{=} \mu^{-2}\{g(t) + (\mu - s/t)^2 h(t)\} + o(1).$$

Tedious calculations show that $g(t) > 0$ and $h(t) > 0$ for all $t > 0$. Thus the mixture rule studied in this section is asymptotically uniformly (in $\mu$) better (by at most an additive constant) than its Robbins–Siegmund analog. This extends a result of Pollak and Yakir [8].

**3. A second example: hypothesis testing and detecting a change of the shape parameter of a Gamma distribution.** As indicated in Section 2, a mixture procedure seems preferable for the normal mean problem. However, there are cases where mixtures are hard to apply, such as distributions that do not admit a conjugate prior, especially when the parameter is multidimensional. In such cases, a Robbins–Siegmund scheme would be of value. In this section, we illustrate this by setting up a power one test and a change-point detection scheme for the shape parameter of a Gamma distribution. The considerations involved are typical of more complex problems.

*A power one test.* Let $X_1, X_2, \ldots$ be i.i.d. Gamma$(\theta, 1)$-distributed random variables, and let $H_0 : \theta = \theta_0$, $H_1 : \theta \neq \theta_0$ where $0 < \theta_0 < \infty$ is fixed. This is an example where there is no "natural" mixture; the Gamma$(\theta, 1)$ family has no conjugate prior.

Following the considerations of Section 2, we need to define a sequence of estimates of $\theta$ that are $F_{n-1}$-measurable. A comparison of estimators will be made in Section 4. In this section we consider a particularly simple yet flexible approach based upon the method of moments. Let $0 \leq s, t < \infty$ and define $\theta_1 = s/t$ $(= \theta_0$ if $s \wedge t = 0)$ and $\theta_{n+1} = (n\bar{X}_n + s)/(n+t)$ for $n \geq 1$. Define a Robbins–Siegmund type of test with

$$\Lambda_n = \prod_{i=1}^{n} [X_i^{\theta_i - \theta_0} \Gamma(\theta_0)/\Gamma(\theta_i)]$$

as its test statistic and

$$\tau_b = \min\{n | n \geq 1, \Lambda_n \geq e^b\}$$

($\tau_b = \infty$ if no such $n$ exists) as its stopping time. For $\theta \neq \theta_0$ let

$$\eta(\theta) = \lim_{b \to \infty} E_\theta \exp\left\{-\left[\sum_{i=1}^{M_b} \log(f_\theta(Z_i)/f_{\theta_0}(Z_i)) - b\right]\right\}$$

where $Z_1, Z_2, \ldots$ are i.i.d. Gamma$(\theta, 1)$-distributed random variables, $f_\theta$ is the Gamma$(\theta, 1)$ density and $M_b = \min\{n | n \geq 1, \sum_{i=1,\ldots,n} \log(f_\theta(Z_i)/f_{\theta_0}(Z_i)) \geq b\}$.

THEOREM 2. *When $\theta_1 = s/t$ $(= \theta_0$ if $s \wedge t = 0)$ and $\theta_{n+1} = (n\bar{X}_n + s)/(n+t)$ for $n \geq 1$, there exists a probability measure $G$ (that depends on $\theta_0, s, t$) on $(0, \infty)$ such that the Robbins–Siegmund test of $H_0: \theta = \theta_0$, $H_1: \theta \neq \theta_0$ has significance level*

(6) $$\alpha = P_{H_0}(\tau_b < \infty) = (1 + o(1)) \times \gamma \times \exp\{-b\},$$

*where $\gamma = \int \eta(y) \, dG(y)$ and $o(1) \to 0$ as $b \to \infty$. The test has power one.*

REMARK. Although the constant $\gamma$ and the measure $G$ do not admit an analytic expression, they can be evaluated easily by Monte Carlo, as will be shown in Section 4. This turns Theorem 2 into a practical tool, as the significance level can be approximately regulated [by letting $b = \log(\gamma/\alpha)$] once $\gamma$ has been evaluated. The choice of $s, t$ influences the ASN when $\theta \neq \theta_0$, as discussed in Section 4.

SKETCH OF PROOF OF THEOREM 2. Under $P_\theta$, clearly $\theta_n \to \theta$ a.s. as $n \to \infty$, so the test has power one, as can be seen easily by the methods of Robbins and Siegmund [9].

The proof of (6) goes along the same lines as that of Theorem 1. Let $Q$ be the measure on $\{X_1, X_2, \ldots\}$ under which the distribution of $X_n$ conditional on $X_1, \ldots, X_{n-1}$ is Gamma$(\theta_n, 1)$. We will prove an analog of Lemma 1. The rest of the proof of (6) is essentially the same as that of Theorem 1, so we omit the details.

LEMMA 1*. *Under the measure $Q$, the sequence $\{\theta_n\}$ converges a.s. to a positive random variable whose distribution is $G$.*

PROOF. By direct calculation, note that under $Q$ the sequence $\{\theta_n\}$ is a martingale with expectation $\theta_1 = s/t$ ($= \theta_0$ if $s \wedge t = 0$). Therefore $E_Q(\exp\{-\theta_{n+1}\}|F_{n-1}) \geq \exp\{-E_Q(\theta_{n+1}|F_{n-1})\} = \exp\{-\theta_n\}$. Thus $\exp\{-\theta_n\}$ is a bounded submartingale under $Q$ and therefore it has an a.s. limit. Consequently, $\{\theta_n\}$ has an a.s. limit. It is left to prove that this limit is a.s. positive and finite. (If the limit were concentrated on 0 and $\infty$, then Theorem 2 would not be of practical value.)

Note that
$$\mathrm{Var}_Q\,\theta_n = \mathrm{Var}_Q\,E_Q(\theta_n|F_{n-2}) + E_Q\,\mathrm{Var}_Q(\theta_n|F_{n-2})$$
$$= \mathrm{Var}_Q\,\theta_{n-1} + E_Q\theta_{n-1}/(t+n-1)^2$$
$$= \mathrm{Var}_Q\,\theta_{n-1} + \theta_1/(t+n-1)^2,$$

so that
$$\mathrm{Var}_Q\,\theta_n = \sum_{i=1}^{n}(\mathrm{Var}_Q\,\theta_i - \mathrm{Var}_Q\,\theta_{i-1}) < \sum_{i=1}^{\infty}[\theta_1/(t+i)^2].$$

Therefore
$$P_Q(\theta_n > x) < \left\{\sum_{i=1}^{\infty}[\theta_1/(t+i)^2]\right\}\Big/(x-\theta_1)^2 \xrightarrow[x\to\infty]{} 0,$$

so that the limiting distribution of $\theta_n$ does not have an atom at $\infty$.

Now consider $\varphi_n(y) = E_Q\exp\{-y\theta_n\}$ for $y > 0$. To show that the limiting distribution of $\theta_n$ does not have an atom at 0, it suffices to show that $[\lim_{n\to\infty}\varphi_n(y)] \to 0$ as $y \to \infty$. Each $\varphi_n(y)$ is decreasing in $y$ and by the bounded convergence theorem is seen to be continuous and to have limit zero at $+\infty$. Denote the inverse function by $\varphi_n^{-1}$. If $0 < \varepsilon < 1$ and the sequence $\{\varphi_n^{-1}(\varepsilon)\}_{n=1,\ldots,\infty}$ is bounded above by $M$ (say), then, for $y > M$, $\lim_{n\to\infty}\varphi_n(y) \leq \lim_{n\to\infty}\varphi_n(\varphi_n^{-1}(\varepsilon)) = \varepsilon$. It therefore suffices to show that, for each $0 < \varepsilon < 1$, $\{\varphi_n^{-1}(\varepsilon)\}$ is bounded above.

Recalling that $E(\exp\{-r \times \mathrm{Gamma}(\theta,1)\}) = (1+r)^{-\theta}$,
$$\varphi_{n+1}(y) = E_Q E_Q\left(\exp\left\{-y\left(s + \sum_{i=1,\ldots,n-1} X_i + X_n\right)\Big/(t+n)\right\}\Big|F_{n-1}\right)$$
$$= E_Q(\exp\{-y\theta_n(t+n-1)/(t+n) - \theta_n\log(1+y/(t+n))\})$$
$$= \varphi_n(y(t+n-1)/(t+n) + \log(1+y/(t+n))).$$

Defining $y_n = \varphi_n^{-1}(\varepsilon)$ and setting $y = y_{n+1}$ in the previous line yields
$$\varphi_n(y_n) = \varepsilon = \varphi_{n+1}(y_{n+1})$$
$$= \varphi_n(y_{n+1}(t+n-1)/(t+n) + \log(1+y_{n+1}/(t+n)))$$

and therefore

$$y_n = y_{n+1}(t+n-1)/(t+n) + \log(1+y_{n+1}/(t+n)).$$

It remains to show that $\{y_n\}$ is bounded above. Letting $\gamma_{n+1} = y_{n+1}/(t+n)$,

(7) $$\gamma_n = \gamma_{n+1} + (\log(1+\gamma_{n+1}))/(t+n-1).$$

Clearly, $\{\gamma_n\}$ is decreasing and must have a limit $q \geq 0$. Hence

$$q - \gamma_2 = \sum_{n=2,\ldots,\infty} (\gamma_{n+1} - \gamma_n) = - \sum_{n=2,\ldots,\infty} (\log(1+\gamma_{n+1}))/(t+n-1) > -\infty.$$

Since $\log(1+\gamma_{n+1}) \to \log(1+q)$, evidently $q = 0$. Thus

$$\gamma_2 = \sum_{n=2,\ldots,\infty} (\log(1+\gamma_{n+1}))/(t+n-1),$$

and since $\log(1+\gamma_{n+1}) > (1/2)\gamma_{n+1}$ for large $n$,

(8) $$\sum_{n=2,\ldots,\infty} \gamma_{n+1}/(t+n-1) < \infty.$$

For all $n$

$$\log(1+\gamma_{n+1}) > \gamma_{n+1} - \gamma_{n+1}^2,$$

and by (7)

$$\gamma_n > \gamma_{n+1} + (\gamma_{n+1} - \gamma_{n+1}^2)/(t+n-1)$$
$$= \gamma_{n+1}(t+n)/(t+n-1) - \gamma_{n+1}^2/(t+n-1).$$

Multiplying by $t+n-1$,

$$y_n > y_{n+1} - \gamma_{n+1}^2 = y_{n+1}(1 - \gamma_{n+1}/(t+n)),$$

so that for $n > k$ sufficiently large the right-hand side is positive and

$$y_k/y_{m+1} = \prod_{n=k,\ldots,m} y_n/y_{n+1} > \prod_{n=k,\ldots,m} (1 - \gamma_{n+1}/(t+n))$$
$$> \prod_{n=k,\ldots,\infty} (1 - \gamma_{n+1}/(t+n)) > 0$$

by (8). Thus $y_{m+1}$ is bounded above, as required. $\square$

THEOREM 3. *Define the Fisher and Kullback–Leibler information numbers*

$$I(\theta) = -E_\theta[\partial^2(\log f_\theta(X))/\partial\theta^2] = d^2(\log\Gamma(\theta))/d\theta^2,$$
$$I(\theta,\phi) = E_\theta \log[f_\theta(x)/f_\phi(x)] = (\theta-\phi)d(\log\Gamma(\theta))/d\theta - \log[\Gamma(\theta)/\Gamma(\phi)].$$

Let $\{\theta_n\}$ be a sequence of $F_{n-1}$-measurable estimators of $\theta$ with asymptotic efficiency $\kappa(\theta)$ in the sense that

$$E_\theta(\theta_n - \theta)^2 = (1 + o(1))/[nI(\theta)\kappa(\theta)] \qquad \text{as } n \to \infty \text{ for all } \theta.$$

Assume that

$$E_\theta(\theta_n - \theta)^4 = O(n^{-2}) \qquad \text{as } n \to \infty$$

and that there exists $c > 0$ depending on $\theta$ such that

$$\sum_{n=1,\ldots,\infty} E_\theta[(\log \theta_n^{-1})^+ \mathbb{1}(\theta_n < c)] < \infty.$$

Let $\tau_b$ be defined by the estimator sequence $\{\theta_n\}$. Then for all $\theta$, as $b \to \infty$

$$E_\theta \tau_b = I(\theta, \theta_0)^{-1}(b + (\log b))/(2\kappa(\theta)) + o(\log b).$$

SKETCH OF PROOF. To make the writing easier, use $N$ as a shorthand for $\tau_b$. Standard estimates show that

$$E_\theta \log \Lambda_N = b + O(1),$$

and letting $\Lambda_n^\theta = f_\theta(X_1, \ldots, X_n)/f_{\theta_0}(X_1, \ldots, X_n)$, it follows using Wald's equation that

$$I(\theta, \theta_0) E_\theta N = b + O(1) + E_\theta(\log \Lambda_N^\theta - \log \Lambda_N).$$

Applying the martingale optional sampling theorem to $\{\log \Lambda_n^\theta - \log \Lambda_n - \sum_{k=1,\ldots,n} I(\theta, \theta_k)\}$, it remains to show that

$$E_\theta\left(\sum_{k=1,\ldots,N} I(\theta, \theta_k)\right) = (\log b)/(2\kappa(\theta)) + o(\log b).$$

Fix $0 < c < \theta$. To verify that there is an $A > 0$ such that for all $\phi$

$$|I(\theta, \phi) - \tfrac{1}{2}(\phi - \theta)^2 I(\theta)| \leq A|\phi - \theta|^3 + (\log \phi^{-1})^+ \mathbb{1}(\phi < c),$$

note the following: the inequality holds for $c \leq \phi \leq \theta + c$ by Taylor expansion of $\log \Gamma(\phi)$ about $\phi = \theta$; it holds for $0 < \phi < c$ since

$$\log \Gamma(\phi) = \log \Gamma(1 + \phi) - \log \phi \leq \text{const} + (\log \phi^{-1})^+;$$

and it holds for $\phi > \theta + c$ since Stirling's approximation yields as $\phi \to \infty$,

$$I(\theta, \phi) \leq \phi \log \phi + O(\phi) \leq O(|\phi - \theta|^3).$$

Thus

$$\left| E_\theta \sum_{k=1,\ldots,N} (I(\theta, \theta_k) - \tfrac{1}{2}(\theta_k - \theta)^2 I(\theta)) \right|$$
$$\leq \sum_{k=1,\ldots,\infty} (A|\theta_k - \theta|^3 + (\log \theta_k^{-1})^+ \mathbb{1}(\theta_k < c)),$$

and since the hypotheses imply that the series converges, it suffices to show that

$$E_\theta\left(\sum_{k=1,\ldots,N}(\theta_k-\theta)^2\right) = (\log b)/[I(\theta)\kappa(\theta)] + o(\log b).$$

Routine modifications of the arguments used to prove Lemmas 13, 14 and 16 of Robbins and Siegmund [10] establish that

$$E_\theta N \sim b/I(\theta,\theta_0) \qquad \text{as } b \to \infty$$

and that by using the definition of $\kappa(\theta)$,

$$E_\theta\left(\sum_{k=1,\ldots,N}(\theta_k-\theta)^2\right) \sim \sum_{k=1,\ldots,[bI(\theta,\theta_0)^{-1}]} E_\theta(\theta_k-\theta)^2 \sim (\log b)/[I(\theta)\kappa(\theta)].$$

$\square$

*Detecting a change.* Now we suppose that when the process being monitored is in control, it yields i.i.d. Gamma-distributed observations, and when the process is out of control the shape parameter changes. Formally stated, an abrupt change may occur at time $\nu$, in which case $X_1, X_2, \ldots, X_{\nu-1}$ are i.i.d. Gamma$(\theta_0, 1)$-distributed random variables and $X_\nu, X_{\nu+1}, \ldots$ are i.i.d. Gamma$(\theta, 1)$-distributed random variables (which are independent of the first $\nu - 1$ observations). The initial shape parameter $\theta_0$ is assumed to be known, but the post-change parameter $\theta$ as well as the changepoint $\nu$ are unknown. We will let $P_\nu$ and $E_\nu$ denote probability and expectation under this scheme, where $\nu = \infty$ denotes no change ever taking place. If the post-change $\theta$ were known, the Shiryayev–Roberts changepoint detection scheme would define the sequence of statistics $R_n^\theta = \sum_{k=1,\ldots,n} f_\theta(X_k,\ldots,X_n)/f_{\theta_0}(X_k,\ldots,X_n)$ and the associated stopping time $N_A^\theta = \min\{n | R_n^\theta \geq A\}$. The sequence $\{R_n^\theta - n\}$ is a $P_\infty$-martingale with zero expectation, a structure used to evaluate the ARL to false alarm of $N_A^\theta$. When $\theta$ is unknown, we propose to estimate it in a way that will preserve the martingale structure. Again the idea is to substitute $F_{n-1}$-measurable estimates for the $\theta$ used in the likelihood ratio of $X_n$.

We present two examples. The first uses a method-of-moments estimator for $\theta$, which leads to simple calculations and a correspondingly simple exposition of the issues involved in proofs. Our second example uses maximum likelihood estimation, which is asymptotically efficient but requires more calculation to apply and more delicate mathematical analysis. We provide a Monte Carlo comparison of the two methods in the next section.

EXAMPLE (*An SRRS procedure based on estimation by the method of moments*). Given $s, t \geq 0$, define

$$\theta_{n,k} = \left( \sum_{i=k,\ldots,n-1} X_i + s \right) \Big/ (n - k + t)$$

for $n \geq k$, where $\theta_{k,k} = \theta_0$ if $s \wedge t = 0$,

(9)
$$\Lambda_{n,k} = \prod_{i=k,\ldots,n} [f_{\theta_{i,k}}(X_i)/f_{\theta_0}(X_i)],$$

$$R_n = \sum_{k=1,\ldots,n} \Lambda_{n,k},$$

$$N_A = \min\{n | R_n \geq A\},$$

where the stopping threshold $A$ is fixed.

Results regarding the operating characteristics of this SRRS procedure are stated in Theorem 4, whose proof is given in the Appendix. Part (iii) of Theorem 4 illustrates the effect of the asymptotic efficiency of the estimation procedure on the delay to detection.

THEOREM 4. *For a Shiryayev–Roberts–Robbins–Siegmund scheme defined by* (9),

(i) $E_\infty N_A \geq A$ *for all* $0 < A < \infty$,
(ii) $\lim_{A \to \infty} E_\infty N_A / A = 1/\gamma$, *where $\gamma$ is the same as in Theorem* 2,
(iii) $\sup_\nu E_{\theta,\nu}(N_A - \nu + 1 | N_A \geq \nu) = \{\log A + \frac{1}{2}(\log \log A)/\kappa(\theta) + o(\log \log A)\}/E_\theta \log(f_\theta(X)/f_{\theta_0}(X))$,

*where* $\kappa(\theta) = 1/(\theta I(\theta))$ *and* $I(\theta)$ *is Fisher information*.

REMARK. Theorem 4 provides a basis for applying this Shiryayev–Roberts–Robbins–Siegmund changepoint detection scheme. If one requires an ARL to false alarm of at least $B$, one can obtain a (conservative) scheme by setting $A = B$, or, after evaluating $\gamma$, a scheme that approximately satisfies the condition by setting $A = B\gamma$. It is possible to obtain asymptotic expressions for the expected delay to detection, but they have constants which have to be evaluated by Monte Carlo separately for each post-change parameter value $\theta$, and since the expressions do not yield good enough approximations for cases of applied interest ($B$ of the order of magnitude $10^2$–$10^3$), we do not present them here.

EXAMPLE (*An efficient estimating sequence*). Theorem 3 implies that better performance can be realized if the estimating sequence is efficient. In this subsection, we apply the (efficient) maximum likelihood estimator sequence instead of the method-of-moments type of sequence studied in the previous two sections. In Theorem 5 we give a formal definition of the procedure and state results regarding its operating characteristics. Proofs appear in the Appendix.

Let $\hat{Q}$ be the measure on $\{X_1, X_2, \ldots\}$ under which $X_1 \sim \text{Gamma}(\theta_0, 1)$ and for $n > 1$ the distribution of $X_n$ conditional on $X_1, \ldots, X_{n-1}$ is $\text{Gamma}(\hat{\theta_n}, 1)$, where $\hat{\theta_n} = $ the solution of $(\sum_{i=1,\ldots,n-1} \log X_i)/(n-1) = E_\theta \log X$, which is easily seen to be the maximum likelihood estimate of $\theta$ based on $X_1, \ldots, X_{n-1}$.

THEOREM 5. (a) *Under the measure $\hat{Q}$, the sequence $\{\hat{\theta_n}\}$ converges a.s. to a positive random variable whose distribution we denote by $G$.*
(b) *For a Shiryayev–Roberts–Robbins–Siegmund scheme defined by*

$$\theta_{n,k} = \text{solution of} \left(\sum_{i=k,\ldots,n-1} \log X_i\right)\Big/(n-k) = E_\theta \log X$$

*for $n \geq k$, where $\theta_{k,k} = \theta_0$,*

$$\Lambda_{n,k} = \prod_{i=k,\ldots,n} [f_{\theta_{i,k}}(X_i)/f_{\theta_0}(X_i)],$$

$$R_n = \sum_{k=1,\ldots,n} \Lambda_{n,k},$$

$$N_A = \min\{n | R_n \geq A\},$$

*the following hold:*

(i) $E_\infty N_A \geq A$ *for all $0 < A < \infty$.*
(ii) $\lim_{A \to \infty} E_\infty N_A / A = 1/\gamma$, *where $\gamma = \int \eta(\theta)\, dG(\theta)$ [$\eta(\cdot)$ is defined in Section 3].*
(iii) $\sup_\nu E_{\theta,\nu}(N_A - \nu + 1 | N_A \geq \nu) = \{\log A + \frac{1}{2}(\log \log A) + o(\log \log A)\}/E_\theta \log(f_\theta(X)/f_{\theta_0}(X)).$

REMARK. For the maximum likelihood procedure, one can retain the flexibility of the method-of-moments produced by introducing the parameters $s, t$ by defining $\theta_{n,k} = $ solution of $(s + \sum_{i=k,\ldots,n-1} \log X_i)/(t + n - k) = E_\theta \log X$. Also, it may make sense to bound the allowable set of $\theta$'s to be bounded away from 0 and from $\infty$, and perhaps also from $\theta_0$. Although this may make the expected delay to detection inefficient for the truncated parameter values, one can argue that they are not of practical interest, whereas their truncation will improve this scheme's performance for the retained set of parameter values.

REMARK. For all SRRS changepoint detection procedures designed for the case that the $P_\infty$-distribution is known, $\sup_\nu E_{\theta,\nu}(N_A - \nu + 1 | N_A \geq \nu)$ is attained at $\nu = 1$. The reason for this is that the $P_{\nu=j}$-behavior of $\{\sum_{k=j,\ldots,n} \Lambda_{n,k}\}_{n=j,\ldots,\infty}$ conditional on $F_{j-1}$ is the same as the $P_{\nu=1}$-behavior of $\{\sum_{k=1,\ldots,n} \Lambda_{n,k}\}_{n=1,\ldots,\infty}$. Let $N_A^{(j)} = \min\{n | n \geq j, \sum_{k=j,\ldots,n} \Lambda_{n,k} \geq A\}$. Clearly, $E_{\nu=1} N_A = E_{\nu=j}(N_A^{(j)} - j + 1 | N_A \geq j)$. However, $R_n \geq \sum_{k=j,\ldots,n} \Lambda_{n,k}$ for all $n \geq j$. Therefore, $N_A^{(j)} \geq N_A$ on $\{N_A \geq j\}$, so that for all $j \geq 1$

$$E_{\nu=1}(N_A - 1 + 1) = E_{\nu=1} N_A = E_{\nu=j}(N_A^{(j)} - j + 1 | N_A \geq j)$$
$$\geq E_{\nu=j}(N_A - j + 1 | N_A \geq j).$$

**4. Monte Carlo.** In this section we present a numerical illustration of the methods proposed in the previous section. We suppose that the prechange distribution is standard exponential, that is, $\theta_0 = 1$. First, we consider the schemes defined by $\theta_1 = 1, \theta_{n+1} = (n\bar{X}_n + t)/(n + t)$ (i.e., $s = t$) for $t = 0, 0.5, 1$.

The first step is to evaluate the constant $\gamma$ (see Theorem 2). For each value of $t$, 5000 replications of $\exp\{-(\log(\Lambda_{\tau_b}) - b)\}$ were run for several ranges of $b$ with $\{X_i\}$ distributed under the measure $Q$ of Lemma 1*. The average of these replications is our estimate of $\gamma$. The results are summarized in Table 1. (See the Appendix for a detailed description of the method of simulation.)

It seems that for $10 < b < 25$ the expectation of $\exp\{-(\log(\Lambda_{\tau_b}) - b)\}$ is nearly constant and hence presumably close to its limiting value. We obtain that $\gamma \approx 0.425, 0.547, 0.606$ for $t = 0, 0.5, 1$, respectively. Next, we ran 10,000 replications to calibrate the Shiryayev–Roberts–Robbins–Siegmund detection scheme to have 500, 750, 1000 as the ARL to false alarm. The results are summarized in Table 2. By Theorem 3, the ratio of $A$ to the ARL to false alarm is asymptotically equal to $\gamma$, and, judging by Table 2, the values of $A$ seem to be large enough for the asymptotics to yield good approximations. In other words, setting $A = \gamma \times$ (desired ARL to false alarm) will achieve the desired ARL to within very few percent. This makes Theorem 3 a practical tool: rather than calibrating $A$ for each ARL separately (which is computationally demanding), it is enough to run a simulation to evaluate $\gamma$ (which takes just a minute or two) and multiply the result by the desired ARL.

TABLE 1
*Monte Carlo estimates of $\gamma$ for three values of $t$, averaged over three intervals*

| $b$-interval | $t = 0$ | | $t = 0.5$ | | $t = 1$ | |
|---|---|---|---|---|---|---|
| $[B_0\ B_1]$ | est. $\gamma$ | s.e. | est. $\gamma$ | s.e. | est. $\gamma$ | s.e. |
| [10 15] | 0.4290 | 0.0044 | 0.5472 | 0.0039 | 0.6065 | 0.0035 |
| [15 20] | 0.4256 | 0.0044 | 0.5502 | 0.0039 | 0.6050 | 0.0036 |
| [20 25] | 0.4215 | 0.0044 | 0.5430 | 0.0040 | 0.6061 | 0.0036 |

TABLE 2
*Levels A (evaluated by Monte Carlo) designed to achieve desired ARLs to false alarm and their relation to $\gamma$, for various values of t*

| $t$ | 0 | | | 0.5 | | | 1 | | |
|---|---|---|---|---|---|---|---|---|---|
| ARL to false alarm | 500 | 750 | 1000 | 500 | 750 | 1000 | 500 | 750 | 1000 |
| $A$ | 221 | 320 | 440 | 275 | 410 | 555 | 309 | 456 | 578 |
| $A/\text{ARL}$ | 0.442 | 0.427 | 0.440 | 0.550 | 0.547 | 0.555 | 0.618 | 0.608 | 0.578 |
| $\gamma$ | | 0.425 | | | 0.547 | | | 0.606 | |

Table 3 presents a comparison of the (maximal) expected delay to detection of three methods-of-moments-based procedures of the kind described in Theorem 4 ($s = t$ and $t = 0, 0.5, 1$) calculated as the average of 10,000 run lengths to detection (when the change is in effect from the start) for each of the post-change parameter values $\theta = 0.35, 0.5, 0.65, 0.8, 1.25, 1.5, 1.75, 2, 2.5, 3$, for ARL to false alarm 1000. The differences are not dramatic, though the choice $t = 1$ seems to give overall performance slightly better than the others.

Also included in Table 3 is a simulation study of the maximum-likelihood-based scheme. The maximum-likelihood-based scheme performs slightly better overall, though it has larger delays to detection when the post-change $\theta$ is less than 1. The calculation of the many maximum likelihood estimates required to perform the SRRS procedure is of course considerably slower than the calculation of the method-of-moments estimators. For each $k$ and $n$ the estimate is obtained by solving numerically for the value of $\theta$ such that $\Gamma'(\theta)/\Gamma(\theta)$ equals the average of $\log(X_k), \ldots, \log(X_n)$.

A central question to be answered is how well do the procedures proposed here compare with simple schemes. For example, since the problem we have been considering is a two-sided problem (the post-change value of $\theta$ may be either larger or smaller than $\theta_0$), a simple alternative method is to choose two values $0 < \theta_1 < \theta_0 < \theta_2 < \infty$, put a prior probability of 50% on $\theta_1$ and on $\theta_2$, and apply the corresponding Shiryayev–Roberts rule; that is, the control statistic is

$$R(n) = \tfrac{1}{2}R_{\theta_1}(n) + \tfrac{1}{2}R_{\theta_2}(n),$$

where $R_{\theta_j}(n)$ are the usual simple Shiryayev–Roberts statistics designed for detecting a change from $\theta_0$ to $\theta_j$ ($j = 1, 2$); that is,

$$R_{\theta_j}(n) = \sum_{k=1,\ldots,n} \Lambda_{k,\theta_j}(n), \qquad j = 1, 2,$$

where

$$\Lambda_{k,\theta_j}(n) = \prod_{i=k}^{n}[X_i^{\theta_i - \theta_0}\Gamma(\theta_0)/\Gamma(\theta_j)], \qquad j = 1, 2,$$

and the stopping time is

$$S_A = \min\{n | R(n) \geq A\}.$$

Based on 10,000 runs, we calibrated $A$ to yield ARL to false alarm 1000 for each of the three pairs $(\theta_1, \theta_2) = (0.8, 1.25), (0.65, 1.5), (0.5, 2)$ (with $\theta_0 = 1$), and ran 10,000 simulations of delay to detection (when the change is in effect from the start) for each of the post-change parameter values $\theta = 0.35, 0.5, 0.65, 0.8, 1.25, 1.5, 1.75, 2, 2.5, 3$, as above. The results are included in Table 3. Not surprisingly, the farther apart $\theta_1$ and $\theta_2$ are, the shorter the expected delay to detection for extreme values of the post-change parameter and the longer the delay for values close to $\theta_0$. Of the three examples chosen, the one with $(\theta_1, \theta_2) = (0.65, 1.5)$ seems most similar to the SRRS ($t = 1$) scheme, which for the values of $\theta$ chosen requires about 15% longer delay for $\theta$ near $\theta_0$ and about 10–30% shorter delay for extreme values of $\theta$.

Finally, in Figure 1 we present histograms of the distribution $G$ of the limit as $n \to \infty$ of $\theta_n$ under the measure $Q$ (for $t = 0, 0.5, 1$). The intervals on the horizontal ($\theta$) axis have width 0.1. This $G$ can be regarded as a "natural" prior on the post-change $\theta$, which could have been used as a mixing measure had mixtures been technically feasible.

It is important to note that we are not trying to make a case for the SRRS procedure as the method of choice for the specific problems considered in

TABLE 3

| Procedure | $t = 0$ | $t = 0.5$ | $t = 1$ | MLE | P-M (0.8, 1.25) | P-M (0.65, 1.5) | P-M (0.5, 2) |
|---|---|---|---|---|---|---|---|
| True $\alpha$ | | | | | | | |
| 0.35 | 10.2 | 9.2 | 9.5 | 9.6 | 15.4 | 10.1 | 8.2 |
| 0.50 | 18.9 | 17.6 | 17.7 | 18.2 | 25.4 | 17.5 | 15.3 |
| 0.65 | 40.2 | 38.0 | 37.2 | 38.6 | 43.7 | 33.6 | 33.4 |
| 0.80 | 112.7 | 104.9 | 101.6 | 108.5 | 94.9 | 94.0 | 122.3 |
| 1.25 | 107.6 | 108.4 | 105.9 | 98.2 | 93.2 | 94.4 | 150.3 |
| 1.50 | 40.8 | 41.6 | 41.1 | 39.4 | 48.0 | 36.0 | 40.1 |
| 1.75 | 23.6 | 24.3 | 24.3 | 23.3 | 34.8 | 23.6 | 20.5 |
| 2.00 | 16.6 | 17.0 | 17.1 | 16.3 | 28.7 | 18.5 | 14.2 |
| 2.50 | 10.2 | 10.6 | 10.8 | 10.3 | 22.4 | 13.8 | 9.6 |
| 3.00 | 7.5 | 7.8 | 8.0 | 7.6 | 19.3 | 11.6 | 7.6 |
| $A$ | 440 | 555 | 578 | 632 | 838 | 700 | 565 |

Monte Carlo: Each cell 10,000 runs; $B = $ ARL to false alarm $= 1000$. (Worst) Average delay for various procedures and various values of (true) $\alpha$.
Estimation: $\theta_{n,k} = (\sum_{i=k, k+1, \ldots, n-1} X_i + t)/(n - k + t)$.
Pair-Mixture: $(\alpha^*, \alpha^{**})$; $P(\alpha = \alpha^*) = P(\alpha = \alpha^{**}) = 1/2$.

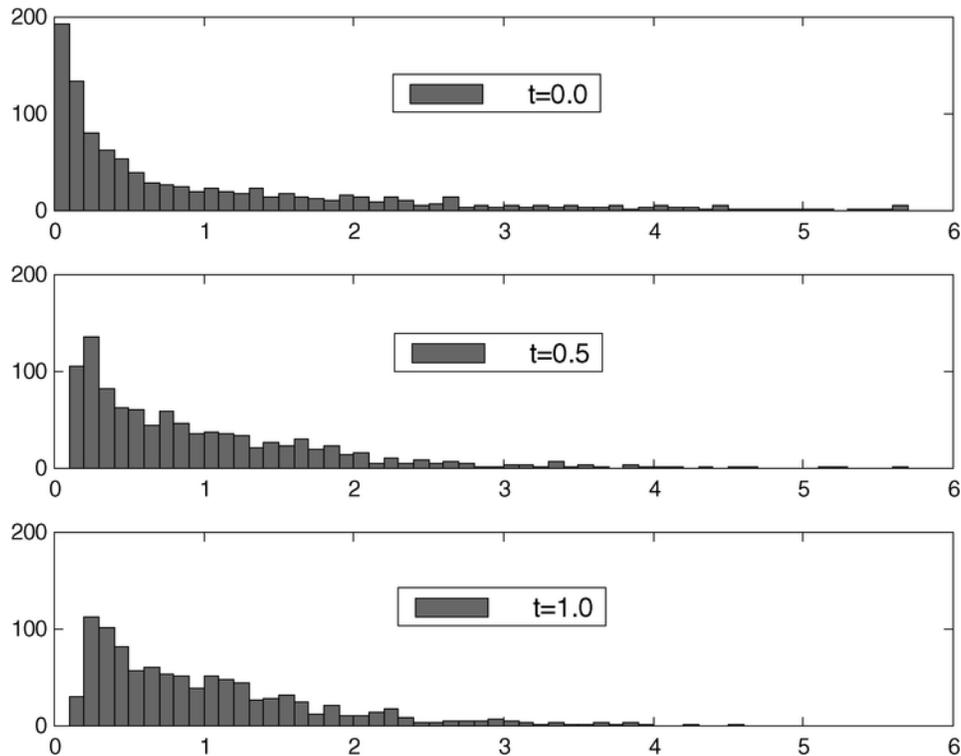

Fig. 1. *Histograms of simulation of G*.

this section and the preceding section. Rather, these sections are meant to introduce the SRRS schemes, show how to apply them, illustrate issues related to proofs of their asymptotics and indicate that they can be effectively used in situations where mixtures are difficult to handle, resulting in at most a slight decrease in efficiency. For single-parameter problems like the Gamma shape parameter, even simpler procedures may show quite acceptable performance (though not asymptotic efficiency), as illustrated by the simulation results for well-chosen mixtures of two simple Shiryayev–Roberts statistics. However, in multiparameter problems such approaches lose their simplicity and transparency, and mixtures are often intractable, in which case the SRRS approach offers worthwhile advantages.

**5. A multiparameter example.** Consider a situation where one monitors simultaneously $m$ independent processes, the $i$th of which yields independent Gamma$(\theta_0 = 1, \beta_i)$-distributed observations when the processes are under control ($\beta_i$ are known), and, either by design or by the nature of the problem, the observations are taken one $m$-vector at a time, the first being observations taken simultaneously from the processes $1, 2, \ldots, m$, the second starting with a second set of observations from the processes $1, \ldots, m$, and so on. A change, if it takes place, may affect some or all of the parameter values, which may be different for the different processes.

For illustration's sake, imagine that the observations are taken one a day, and, when everything is under control, the distribution of an observation is exponential with a parameter that depends on the day of the week in which the observation was made. After standardization, all of the observations have a Gamma$(\theta_0 = 1, 1)$ distribution (pre-change). If there is a changepoint, then subsequent observations are assumed to have a Gamma$(\theta, 1)$ distribution, where the post-change shape parameter depends on the day of the week. In other words, changes may be in the $\theta$ value for one of the days, for some of them, or even for all of them, and the post-change parameter values may differ for different days. (The observations are all assumed to be independent.) We assume that the observations are obtained weekly, and a change may take place only between weeks. Here $m = 7$, and the observations are vectors $\mathbf{X_i}$, where $X_{ij}$ denotes the observation on the $j$th day of the $i$th week. For a method-of-moments approach, define

$$\theta_{n,k,j} = \left( \sum_{i=k,\ldots,n-1} X_{ij} + s \right) \Big/ (n - k + t)$$

for $n \geq k$, where $\theta_{k,k,j} = \theta_0 = 1$ if $s \wedge t = 0$,

$$\Lambda_{n,k} = \prod_{i=k,\ldots,n; j=1,\ldots,7} [f_{\theta_{i,k,j}}(X_{ij})/f_{\theta_0}(X_{ij})],$$

$$R_n = \sum_{k=1,\ldots,n} \Lambda_{n,k},$$

$$N_A = \min\{n | R_n \geq A\}.$$

For a maximum-likelihood approach, take $\theta_{n,k,j} = $ solution of $E_\theta \log X = (s + \sum_{i=k,\ldots,n-1} \log X_{ij})/(n - k + t)$ for $n > k$, and $\theta_{k,k,j} = \theta_0 = 1$.

Analogs of Theorems 4 and 5 are valid, the only change being that $\gamma$ has to be recalculated (by Monte Carlo, in a manner analogous to Section 4). The application of the schemes is straightforward, requiring a short computer program.

In this example, even a discrete mixture is impractical. The simplest reasonable choice would put a prior of 1/3, 1/3, 1/3 on $\theta_0 = 1$ and some $\theta_1$ and $\theta_2$, independently for each of the seven days of the week, leading to a discrete mixture of $3^7 = 2187$ points [deleting, perhaps, the point $(\theta_0, \theta_0, \ldots, \theta_0)$]. (One needs at least three $\theta$'s to allow for the fact that there may be a change in only a subset of the parameters.) Schemes that put weight on more than three $\theta$'s are even more cumbersome. Schemes that reduce the number of points will be inefficient for detecting certain constellations of change. On the other hand, the SRRS scheme is intuitive and fairly easy to implement. Furthermore, it has the flexibility of easily accommodating prior knowledge of the region where post-change parameters may be. For example, if the only possibility of a change is for the shape parameter to increase, the estimator can be restricted to $\theta > \theta_0$. Or, if for some reason it is clear that the post-change shape parameter will be larger on weekdays than on weekends, then the estimators may be modified to reflect this. This may be a reason to consider SRRS schemes even in problems where in the nonrestricted settings mixtures are feasible; if restrictions are added, integration with respect to a conjugate prior may prove to be much less tractable than in an unrestricted context.

**6. Remarks.** 1. Intuitively, one would expect the Robbins–Siegmund type of rule to do somewhat worse than its mixture counterpart. For example, in the Gamma shape parameter problem, if one takes $s = 0$, then the parameter value used for the first likelihood ratio equals the initial pre-change parameter value, so that the (estimated) likelihood ratio of the first observation is unity no matter what the value of the first observation is. In other words, one "loses" an observation, something which does not occur when employing the mixture analog. The decision whether or not to stop sampling at the $n$th stage is based on a sufficient statistic under the mixture rule but not under the Robbins–Siegmund rule. Nonetheless, the latter method seems to perform nearly as well as the mixture rule, as indicated by Theorem 2 and the discussion at the end of Section 2. We ran a simulation to compare the SRRS and its mixture analog for detecting a change of a normal mean. (Here mixtures are easy to implement; we just wanted to see how the methods compare in a "standard" context.) We assumed the variance of the observations to be 1 and the pre-change mean to be 0. By numerical calculations of the variance of $G_{s,t}$ of Section 2, $G_{s=0,t=0.42626} = N(0,1)$. We constructed the SRRS scheme (with $s = 0$, $t = 0.42626$) in the same manner as is done in (9) for the Gamma example and compared it by simulation to the mixture rule with a $N(0,1)$ mixing measure. The results are recorded in Table 4. (The case $\mu = 0$ gives the simulated ARLs to false alarm. For $\mu > 0$ $\nu$ was taken to be 1.)

Table 4 indicates that the time to detection is, as to be expected, somewhat longer for the SRRS rule, but the difference is not very great. For $\mu = 0.25$ the difference is insignificant, and for the other values of $\mu$ that were investigated, over the range $\mu = 0.5$ to $\mu = 3$ there is a remarkably consistent pattern: the ARLs of the SRRS rule are about 0.4 or 0.5 larger than those of the mixture rule. (The ARLs to false alarm for the two rules are roughly equal for each $A$ in the range investigated, with the mixture rule having an ARL about 1–2% larger than the SRRS. This difference is small, and adjusting for it hardly changes the picture. The overshoot effect is $E_\infty N_A / A \approx 1/\gamma \approx 1.5$.)

2. Although we conjecture that the SRRS scheme is never better than its mixture analog, the following example indicates that it is in some cases no worse.

Let $X_i \sim \text{Binomial}(1,p)$ i.i.d.; $H_0: p = p_0$, $H_1: p \neq p_0$. Suppose $p_1 = s/t$, $p_{n+1} = (s + \sum_{i=1,\ldots,n} X_i)/(t+n)$ where $0 < s < t < \infty$ are constants. Note that the behavior of the sequence $X_1, X_2, \ldots$ [with the conditional distribution of $X_n$ given the past being $\text{Binomial}(1, p_n)$] is identical to that of the sequence $X_1, X_2, \ldots$, where, conditioned on $p$, the $X_i$ are i.i.d. $\text{Binomial}(1,p)$ variables and there is a $\text{Beta}(s, t-s)$ prior on $p$. In this setting clearly $p_n \to p$ a.s. as $n \to \infty$. Hence $G$ is $\text{Beta}(s, t-s)$ (the same as the prior on $p$). In this example SRRS is identical to its mixture counterpart.

3. When there is a suitable invariance structure, the Robbins–Siegmund technique can be applied also when the baseline is unknown. To illustrate this, consider again the change of normal mean problem as in Remark 1, but suppose that the initial baseline (the pre-change mean) is unknown. Invariance considerations would base changepoint detection on the sequence $\{Y_i\} = \{X_i - X_1\}$ instead of the original sequence $\{X_i\}$ (see [7]). The unknown post-change parameter $EY_n$ can be estimated by $(Y_k + \cdots + Y_{n-1})/(n-k)$.

4. Usually there is no obvious natural prior for a mixture rule, whereas there are natural estimates. At least in theory, in such cases the estimates can be regarded as inducing a natural prior. For instance, in the example treated in Section 2, if $\bar{X}$ is considered to be a natural estimate of the mean,

TABLE 4
*Simulated ARLs for detecting a change of a normal mean, 40,000 runs*

| $A$ | $\mu = 0$ | $\mu = 0.25$ | $\mu = 0.5$ | $\mu = 0.75$ | $\mu = 1$ | $\mu = 1.5$ | $\mu = 2$ | $\mu = 3$ |
|---|---|---|---|---|---|---|---|---|
| 400 | 587 | 104.8 | 38.5 | 20.9 | 13.57 | 7.55 | 5.11 | 3.18 |
|  | 599 | 104.7 | 38.1 | 20.5 | 13.13 | 7.11 | 4.68 | 2.73 |
| 450 | 661 | 109.2 | 39.6 | 21.3 | 13.82 | 7.66 | 5.17 | 3.21 |
|  | 673 | 109.0 | 39.1 | 20.9 | 13.38 | 7.22 | 4.74 | 2.76 |
| 500 | 739 | 113.0 | 40.5 | 21.7 | 14.05 | 7.76 | 5.23 | 3.24 |
|  | 748 | 112.9 | 40.1 | 21.3 | 13.60 | 7.32 | 4.80 | 2.78 |
| 550 | 813 | 116.6 | 41.3 | 22.1 | 14.25 | 7.84 | 5.28 | 3.26 |
|  | 823 | 116.3 | 40.9 | 21.7 | 13.80 | 7.41 | 4.85 | 2.81 |
| 600 | 886 | 119.7 | 42.1 | 22.4 | 14.44 | 7.93 | 5.33 | 3.29 |
|  | 900 | 119.7 | 41.6 | 22.0 | 13.97 | 7.49 | 4.90 | 2.83 |
| 650 | 961 | 122.8 | 42.8 | 22.7 | 14.62 | 8.00 | 5.37 | 3.31 |
|  | 973 | 122.7 | 42.3 | 22.3 | 14.14 | 7.57 | 4.94 | 2.85 |
| 700 | 1037 | 125.5 | 43.4 | 23.0 | 14.77 | 8.07 | 5.41 | 3.32 |
|  | 1052 | 125.4 | 43.0 | 22.6 | 14.30 | 7.64 | 4.98 | 2.87 |
| s.e. | 0.43 | 0.35 | 0.11 | 0.05 | 0.03 | 0.014 | 0.008 | 0.004 |

(For each $A$, first row = SRRS, second row = mixture.)

then $s = t = 0$ and $G_{s=0,t=0} = N(0, \pi^2/6)$. So one can argue that a natural mixture rule is one with a $N(0, \pi^2/6)$ prior on $\mu$.

5. For a practical application, it is not imperative to prove an analog of Lemma 1 (though its validity is an ingredient in ensuring asymptotic optimality). Heuristically, the overshoot correction can be expected to be nearly a constant function of $A$ once $A$ is reasonably large, and the constant can be estimated by the simulation methods discussed in Section 4.

6. Another approach to evaluation of $\gamma$, the limit value of the overshoot correction, has been proposed by Yakir and Pollak [19]. That method has the potential to allow error estimates, but in the problem of detecting a change in $\theta$, the shape parameter of the Gamma distribution, the error estimates proved to be difficult to apply.

7. Deleting the last observation from use in the estimation process preserves the martingale structure of the Shiryayev–Roberts statistic. Initially, it seems unnatural to exclude the last observation: after all, this seems to entail a slight loss of efficiency and foregoing sufficiency. The following example shows that there is more involved than mere preservation of a mathematical (martingale) property: at least in the hypothesis testing case, inclusion of the last observation in the estimation sequence can wreak havoc with the level of significance.

As in Section 2, consider a power one test of $H_0 : X_i \sim N(0, 1)$ versus $H_1 : X_i \sim N(\mu, 1)$, where $X_i$ are i.i.d. random variables and $0 \neq \mu \in (-\infty, \infty)$. Here the likelihood ratio is $\Lambda_n(\mu) = \exp(\mu \sum_{i=1,\ldots,n} X_i - n\mu^2/2)$. If one substitutes the maximum likelihood estimator $\sum_{i=1,\ldots,n} X_i/n$ for $\mu$, then the stopping rule based on the estimated likelihood ratio becomes $N_A = \min\{n | \Lambda_n(\sum_{i=1,\ldots,n} X_i/n) \geq A\} = \min\{n | |\sum_{i=1,\ldots,n} X_i| \geq [2n \log(A)]^{1/2}\}$. The law of the iterated logarithm implies that $P_{H_0}(N_A < \infty) = 1$ regardless of the value of $A$, so that one loses the capability of having a nontrivial probability bound on the level of significance.

One stands to lose even if one uses the $n$th maximum likelihood estimator for the $n$th likelihood ratio only—that is, write $\Lambda_n = \exp\{\sum_{i=1,\ldots,n}(\mu_i X_i - \mu_i^2/2)\}$ with $\mu_i = \sum_{j=1,\ldots,i} X_j/i$, and define $N_A = \min\{n | \Lambda_n \geq A\}$. One can show that here, too, $P_{H_0}(N_A < \infty) = 1$ regardless of the value of $A$. (Sketch of proof: Show that $E_{H_0} \log \Lambda_n = \frac{1}{2}(\log n)(1 + o(1))$ and $\text{Var}_{H_0} \Lambda_n = \frac{1}{4}(\log n)^2(1 + o(1))$ as $n \to \infty$. Argue that for $0 < \varepsilon < 1$, asymptotically $P_{H_0}\{\log \Lambda_n > \varepsilon \times \frac{1}{2} \log n\} > \delta$ for some $\delta > 0$. Then break the time axis into intervals $[1, n_1], [n_1 + 1, n_2], [n_2 + 1, n_3], \ldots$ large enough so that $\log \Lambda_{n_i}$ are "almost" independent, and conclude that for any fixed $A$, $P_{H_0}\{\Lambda_n \geq A \text{ for some } 1 \leq n < \infty\} = 1$.)

This phenomenon is not as marked in the changepoint detection context. See [15].

8. In the multiparameter case, our methods are more flexible than indicated. For example, reconsider Section 5. Our methods can be designed for the case that observations are taken on a daily basis, and the change may occur between days rather than only between weeks.

Let the observations be labeled $X_i$, where $i$ is number of days since the onset of changepoint detection, and define $\theta_{n,k} = (s + X_{n-7} + X_{n-14} + \cdots + X_{n-7r})/(t+r)$ [or define $\theta_{n,k,j} =$ solution of $E_\theta \log X = (s + \log X_{n-7} + \log X_{n-14} + \cdots + \log X_{n-7r})/(t+r)$] where $r = r(n,k) = \lfloor (n-k)/7 \rfloor$ and $\Lambda_{n,k}, R_n$ and $N_A$ are as in (9).

Here, too, analogs of Theorems 4 and 5 are valid, with $\gamma$ having to be recalculated (by Monte Carlo, in a manner analogous to Section 4). The application of the schemes is straightforward, requiring a short computer program.

9. The argument used to prove Theorem 5(a) can serve as a model for dealing with many similar problems. The essential ingredients are that the estimators $\theta_{n+1}^*$ satisfy an equation of the form

$$E_{\theta_{n+1}^*} T(X) = [T(X_1) + \cdots + T(X_n)]/n$$

and that an analog of (15) holds, that is,

$$\operatorname{Var}_\theta T(X) \le a + b(E_\theta T(X))^2 \qquad \text{for all } \theta.$$

**7. Conclusion.** We propose that the SRRS scheme is an efficient detection scheme, and should be useful wherever mixture rules are desired but hard to implement. The construction and application of an SRRS rule is simple: all one needs is a sequence of (preferably efficient) estimators for the post-change parameter based on the first $n-1$ observations. Each such estimator will be used to construct an estimated likelihood ratio of the $n$th observation. The likelihood ratios are used to construct a Shiryayev–Roberts statistic, as done in Section 3. In order to achieve an ARL to false alarm $B$, a conservative rule will stop and declare a change to be in effect when the Shiryayev–Roberts statistic first exceeds $B$. A rule that attains $B$ approximately as its ARL to false alarm will stop and declare a change to be in effect when the Shiryayev–Roberts statistic first exceeds $A = B\gamma$. The constant $\gamma$ has to be evaluated (usually) by simulation of tests of hypotheses as in Section 4, but this is the only simulation required, and it takes very little computer time.

## APPENDIX

### A.1. Sketch of proof of Theorem 4.

SKETCH OF PROOF OF THEOREM 4(i). Note that $\{R_n - n\}$ is a $P_\infty$-martingale with zero expectation, and by the optional sampling theorem $E_\infty(R_{N_A} - N_A) = 0$. Since by definition $R_{N_A} \ge A$, this implies that $E_\infty N_A = E_\infty(R_{N_A}) \ge A$, which establishes (i).

SKETCH OF PROOF OF THEOREM 4(ii). We follow along the lines of the proof of Yakir [18], Theorem 3 (and Theorem 1). Before introducing notation in (11) below, we sketch the ideas of the proof.

Break up the time axis (the positive integers) into pieces of size $m$, and show that the $P_\infty$-distribution of $N_A$ can be approximated by using the distribution of the first block (of $m$ observations) where stopping occurs. More precisely, given an integer $j$, the idea is to define (where $A \setminus B$ denotes $A \cap B^c$)

$$\bar{S}_{j,m} = \{jm \leq N_A\}, \qquad S_{j,m} = \bar{S}_{j,m} \setminus \bar{S}_{j+1,m},$$

and to show that

(10) $$(1-\varepsilon)\gamma m/A \leq P_\infty(S_{j,m}|\bar{S}_{j,m}) \leq (1+\varepsilon)\gamma m/A.$$

This enables approximation of $N_A$ by

$$m \times \{\text{a Geometric}(p = \gamma m/A)\text{-distributed random variable}\},$$

from which $E_\infty N_A \approx m/(\gamma m/A) = A/\gamma$ follows.

In order to carry out this program, it turns out that one needs

$$\log A \ll m \ll A.$$

The key ingredient for proving (10) is a measure transformation that will be shown to yield

$$P_\infty(S_{0,m}|\bar{S}_{0,m}) = P_\infty(N_A \leq m)$$
$$= \sum_{k=1,\ldots,m} E_k(\exp\{-(\log R_{N_A} - \log A)\}; \{k \leq N_A \leq m\})/A.$$

Since $\log A \ll m$, for "most" $k$'s $P_k\{k \leq N_A \leq m\} \approx 1$. Also, a renewal-theoretic argument will show that the asymptotic ($A \to \infty$) behavior (under $P_k$) of $(\log R_{N_A} - \log A)$ is the same as that of the log-likelihood ratio statistic in the context of the power one test. Therefore, $P_\infty(S_{0,m}|\bar{S}_{0,m}) \approx m\gamma/A$. The argument is extended to general $P_\infty(S_{j,m}|\bar{S}_{j,m})$ by induction on $j$.

In order to make the analogy to Yakir [18] more transparent, note that the Gamma$(\theta, 1)$ family can be transformed into an exponential family with canonical form: if $X \sim \text{Gamma}(\theta, 1)$, then a reparameterization $y = \theta - \theta_0$ and an appropriate affine transformation $X^*$ of $\log X$ yield a family of probability measures of $X^*$ with densities

$$f_y(x) = \exp\{yx - \psi(y)\}f_0(x), \qquad y \in (-\theta_0, \infty),$$

where $\psi(0) = \psi'(0) = 0$. With this notation, the estimator $y(n, k)$ (for the parameter of the $n$th observation under the putative $\nu = k$) dictated by (9)

is
$$y(n,k) = \left(\sum_{i=k,\ldots,n-1} X_i + s\right)\Big/(n-k+t) - \theta_0$$
$$\text{for } n \geq k, \text{ where } y_{k,k} = 0 \text{ if } s \wedge t = 0.$$

Roughly emulating Yakir's notation, let

$$Z_i^y = yX_i^* - \psi(y),$$

$$dP_k^{y(n,k)}/dP_\infty = \exp\left\{\sum_{i=k,\ldots,n} Z_i^{y(i,k)}\right\},$$

$$R_n = \sum_{k=1,\ldots,n} \exp\left\{\sum_{i=k,\ldots,n} Z_i^{y(i,k)}\right\} = \sum_{k=1,\ldots,n} dP_k^{y(n,k)}/dP_\infty,$$

$$N_A = \min\{n | R_n \geq A\},$$

$$a = \log A,$$

$$M(A) = \min\left\{n \Big| \sum_{i=1,\ldots,n} Z_i^{y(i,1)} \geq \log A\right\} = \tau_A,$$

(11) $\quad Q_k =$ measure analogous to the $Q$-measure of the proof of Theorem 1, appropriate to the Gamma$(\theta, 1)$ context dictated by (9), applied to $X_k^*, X_{k+1}^*, \ldots,$

$H =$ asymptotic distribution (under the measure $Q_{k=1}$) of the overshoot
$$\sum_{i=1,\ldots,\tau_A} Z_i^{y(i,1)} - a,$$

$$\gamma = \int \exp\{-x\} dH(x).$$

We obtain an analog of Yakir's [18] Lemma 1:

LEMMA Y1. *Let $m = m(A)$ be a sequence that satisfies*
$$A/m \to \infty \quad \text{and} \quad (\log A)/m \to 0 \quad \text{as } A \to \infty.$$
*Then*

(12) $\quad P_\infty(N_A \leq m)/P_\infty(M(A/m) \leq m) \to 1 \quad \text{as } A \to \infty.$

SKETCH OF PROOF. For any stopping time $N$,
$$P_\infty(N_A \leq m) = \sum_{j=1,\ldots,m} \int_{\{N=j\}} dP_\infty$$

$$= \sum_{j=1,\ldots,m} \sum_{k=1,\ldots,j} \int_{\{N=j\}} (dP_k^{y(j,k)}/R_j)$$

$$= \sum_{k=1,\ldots,m} \int_{\{k\leq N\leq m\}} (dP_k^{y(N,k)}/R_N).$$

Let $r(n) = \log R_n$. Now

$$(13) \quad \int_{\{k\leq N_A\leq m\}} (dP_k^{y(N,k)}/R_N) = \int_{\{k\leq N_A\leq m\}} \exp\{-(r(N_A) - a)\} \, dQ_k/A.$$

By an argument analogous to the proof of Theorem 1, the denominator in (12) can be shown to be

$$P_\infty(M(A/m) \leq m) = (\gamma m/A)(1 + o(1)),$$

since $m$ grows faster than $\log A$. Therefore, it will suffice to show that, for most of the $k$'s, the value of the integral on the right-hand side of (13) is approximately $\gamma$. Note that

$$R_{k-1+n} = e^{\sum_{i=k,\ldots,k-1+n} Z_i^{y(i,k)}}$$

$$\times \left[ \sum_{j=1,\ldots,k-1} e^{\sum_{i=j,\ldots,k-1} Z_i^{y(i,j)}} e^{-\sum_{i=k,\ldots,k-1+n}(Z_i^{y(i,k)} - Z_i^{y(i,j)})} + 1 \right.$$

$$\left. + \sum_{j=k+1,\ldots,k-1+n} e^{-\sum_{i=k,\ldots,j-1} Z_i^{y(i,k)}} e^{\sum_{i=j,\ldots,k-1+n}(Z_i^{y(i,j)} - Z_i^{y(i,k)})} \right],$$

so that

$$r(k-1+n) - \sum_{i=k,\ldots,k-1+n} Z_i^{y(i,k)} \stackrel{\text{def}}{=} \log[W_0(k,n) + 1 + W_1(k,n)].$$

Observe that for $i > c > b$

$$y(i,c) - y(i,b)$$

$$= \left( \sum_{u=c,\ldots,i-1} X_u + s \right) \Big/ (i-c+t) - \left( \sum_{u=b,\ldots,i-1} X_u + s \right) \Big/ (i-b+t)$$

$$= \left[ (c-b)s + (c-b) \sum_{u=c,\ldots,i-1} X_u - (i-c+t) \sum_{u=b,\ldots,c-1} X_u \right]$$

$$\times [(i-c+t)(i-b+t)]^{-1},$$

and argue that $W_0(k,n) \to W_0(k,\infty)$ and $W_1(k,n) \to W_1(k,\infty)$ a.s. $(Q_k)$ as $n \to \infty$, where both $W_0(k,\infty)$ and $W_1(k,\infty)$ are a.s. $(Q_k)$ finite. Also note that for $r > 0, u_1 > 0, u_2 > 0$,

$$|\log(r+1+u_1) - \log(r+1+u_2)| \leq |\log(1+u_1) - \log(1+u_2)|.$$

These relations imply that Theorem A.7 of Siegmund [14] applies, uniformly in $k$ and in the value of $R_{k-1}$. So (nonlinear renewal theory implies that) the overshoot $r(N_A) - a$, given the value of $R_{k-1}$, has the same asymptotic distribution as the overshoot of the random walk $\sum_{i=1,\ldots,M(A)} Z_i^{y(i,1)}$.

An argument identical to that of Yakir ([18], first half of page 276, verbatim), after replacing Yakir's $R(j,y)$, $r(j,y)$, $N(A,y)$, $dP_k^y$ and $\gamma$ by $R_j$, $r(j)$, $N_A$, $dP_k^{y(N,k)}$ and $\gamma$, completes the proof of Lemma Y1. With the same replacements, the rest of the argument of Yakir (verbatim, from the middle of page 276 until the end of Section 2 on page 278) accounts for our Theorem 4(ii).

SKETCH OF PROOF OF THEOREM 4(iii). One has to check that the conditions of Theorem 3 are satisfied. It is straightforward to check that $E_\theta(\theta_n - \theta)^4 = O(1/n^2)$. As for the other condition, take $c = \theta/2$, write $m = n - 1$ and note that

$$E_\theta(\log^+(\theta_n^{-1}))\mathbb{1}(\theta_n < c)$$
$$\leq \log(n - 1 + t) P_\theta(\theta_n < \theta/2)$$
$$- E_\theta\left(\log\left(\sum_{i=1,\ldots,m} X_i + S\right)\right)\mathbb{1}(\theta_n < \theta/2)$$
$$\times \mathbb{1}\left(\sum_{i=1,\ldots,n-1} X_i + S \leq n - 1 + t\right)\mathbb{1}\left(\sum_{i=1,\ldots,m} X_i + S \leq 1\right),$$

and, since $|\log x| < 1/x$ for $0 < x < 1$,

$$-E_\theta\left(\log\left(\sum_{i=1,\ldots,m} X_i + S\right)\right)\mathbb{1}(\theta_n < \theta/2)$$
$$\times \mathbb{1}\left(\sum_{i=1,\ldots,n-1} X_i + S \leq n - 1 + t\right)\mathbb{1}\left(\sum_{i=1,\ldots,m} X_i + S \leq 1\right)$$
$$\leq \int (1/x)(x^{\theta m - 1}/\Gamma(\theta m))e^{-x}\,dx$$
$$= (1/(\theta m))P(\text{Gamma}(\theta m - 1, 1) < 1).$$

Apply standard manipulations of the Gamma distribution to obtain

$$\sum_{n=1,\ldots,m} E_\theta(\log^+(\theta_n^{-1}))\mathbb{1}(\theta_n < \theta/2) < \infty.$$

**A.2. Sketch of proof of Theorem 5.** Once (a) is proved, the rest follows along the same lines of the proofs in the cases considered in the former sections.

PROOF OF THEOREM 5(a). Let $Y_n = \sum_{i=1,\ldots,n}(\log X_i)/n$ and note that

$$E_{Q^{\hat{}}}(\log X_{n+1}|F_n) = E_{\hat{\theta}_{n+1}}\log X = Y_n.$$

Because $\{Y_n\}$ is a $Q^{\hat{}}$-martingale, general martingale considerations imply that $\mathrm{Var}_{Q^{\hat{}}} Y_n$ increases in $n$ [since $0 \leq E_{Q^{\hat{}}}E_{Q^{\hat{}}}((Y_{n+1}-Y_n)^2|F_n) = \mathrm{Var}_{Q^{\hat{}}} Y_{n+1} - \mathrm{Var}_{Q^{\hat{}}} Y_n$]. Since

$$E_{Q^{\hat{}}}[(\log X_{n+1} - Y_n)^2|F_n] = \mathrm{Var}_{\hat{\theta}_{n+1}} \log X,$$

by writing

$$Y_{n+1} = Y_n + (\log X_{n+1} - Y_n)/(n+1),$$

squaring both sides and taking conditional expectations one obtains

(14)  $$E_{Q^{\hat{}}}(Y_{n+1}^2|F_n) = Y_n^2 + [\mathrm{Var}_{\hat{\theta}_{n+1}} \log X]/(n+1)^2.$$

For $X \sim \mathrm{Gamma}(\theta, 1)$, one obtains by standard considerations that

$$\theta E_\theta \log X \to -1 \quad \text{and} \quad \theta^2 \mathrm{Var}_\theta \log X \to 1 \quad \text{as } \theta \to 0,$$

so there is a $\theta^* > 0$ such that $\mathrm{Var}_\theta \log X \leq 2(E_\theta \log X)^2$ for $\theta \leq \theta^*$. Also, since $\theta \mathrm{Var}_\theta \log X \to 1$ as $\theta \to \infty$ and $\mathrm{Var}_\theta \log X$ is continuous in $\theta$, there is a constant $c$ such that $\mathrm{Var}_\theta \log X \leq c$ for $\theta > \theta^*$. Thus

(15)  $$\mathrm{Var}_\theta \log X \leq c + 2(E_\theta \log X)^2 \quad \text{for all } \theta.$$

Combining (14) and (15) and taking expectations, one obtains

$$E_{Q^{\hat{}}}Y_{n+1}^2 \leq E_{Q^{\hat{}}}Y_n^2 + E_{Q^{\hat{}}}[c + 2(E_{\hat{\theta}_{n+1}}\log X)^2]/(n+1)^2$$
$$= E_{Q^{\hat{}}}Y_n^2 + E_{Q^{\hat{}}}[c + 2E_{Q^{\hat{}}}Y_n^2]/(n+1)^2$$
$$= (1 + 2/(n+1)^2)E_{Q^{\hat{}}}Y_n^2 + c/(n+1)^2.$$

Since $\mathrm{Var}_{Q^{\hat{}}} Y_n$ was shown to be an increasing sequence,

$$E_{Q^{\hat{}}}Y_n^2 \geq \mathrm{Var}_{Q^{\hat{}}} Y_n \geq \mathrm{Var}_{Q^{\hat{}}} Y_1 \geq c/\delta$$

for some $\delta > 0$, and hence

$$E_{Q^{\hat{}}}Y_{n+1}^2 \leq [1 + (2+\delta)/(n+1)^2]E_{Q^{\hat{}}}Y_n^2.$$

This shows that $\{E_{Q^{\hat{}}}Y_{n+1}^2\}$ is bounded, since the infinite product $\prod_{n=1,2,\ldots}[1 + (2+\delta)/(n+1)^2]$ converges. Hence $\{\mathrm{Var}_{Q^{\hat{}}} Y_n\}$ is bounded (and, being monotone, is convergent). It follows that the martingale $\{Y_n\}$ has an a.s. $(Q^{\hat{}})$ finite limit, and consequently $\{\hat{\theta}_{n+1}\}$ has a finite positive limit.

SKETCH OF PROOF OF THEOREM 5(b)(i). As in previous cases, $\{R_n - n\}$ is a $P_\infty$-martingale with zero expectation, and Theorem 5(b)(i) is established by the optional sampling theorem. To see that the conditions of this theorem are met, it suffices to show that $N_A$ is bounded from above by const $\times$ a geometrically distributed random variable. Note that

$$\sum_{k=1}^{n} \Lambda_{n,k} > \Lambda_{n,n-1} = X_n^{\theta_{n,n-1} - \theta_0} \Gamma(\theta_0)/\Gamma(\theta_{n,n-1}) \stackrel{\text{def}}{=} \xi_n.$$

Since $\theta_{n,n-1}$ depends only on $X_{n-1}$, clearly $\xi_n$ depends only on $X_{n-1}$ and $X_n$. Therefore, $\xi_2, \xi_4, \xi_6, \ldots$ are i.i.d. random variables under $P_\infty$, and thus $N_A \leq \min\{n | n = 2m, \xi_n \geq A\}$, which is $2 \times$ a geometrically distributed random variable.

PROOF OF THEOREM 5(b)(ii). The proof follows along the same lines as that of Theorem 4(ii) and is therefore omitted.

SKETCH OF PROOF OF THEOREM 5(b)(iii). It suffices to show that the conditions of Theorem 3 are met.

Let $\zeta(y) = E_y \log X = [d\Gamma(y)/dy]/\Gamma(y)$, and let $\zeta^{-1}$ denote the inverse function of $\zeta$. Observe that $\zeta(y) = \log y + O(1)$ as $y \to \infty$, that $\zeta'(y) = d\zeta(y)/dy = \text{Var}_y \log X = (1 + o(1))/y$ as $y \to \infty$ and that $\text{Var}_y \log X$ is a decreasing function of $y$. Therefore, since $d\zeta^{-1}(t)/dt = 1/\zeta'(\zeta^{-1}(t))$, there exists a finite constant $c_1 > 0$ such that for $t \geq E_{\theta_0} \log X$

$$\begin{aligned}
0 &\leq \zeta^{-1}(t) - \zeta^{-1}(E_{\theta_0} \log X) \\
&\leq (t - E_{\theta_0} \log X) \, d\zeta^{-1}(t)/dt \\
&\leq (t - E_{\theta_0} \log X) c_1 e^t.
\end{aligned} \quad (16)$$

For $p > 0$, let $\delta > \theta_0^p$. Since $\hat{\theta}_{n+1} = \zeta^{-1}(\overline{\log X})$ (where $\overline{\log X} = \sum_{i=1}^{n} \log X_i/n$),

$$E_{\theta_0}|\hat{\theta}_{n+1} - \theta_0|^p$$
$$(17) \qquad = \int_0^\delta P_{\theta_0}(|\hat{\theta}_{n+1} - \theta_0|^p > t) \, dt + \int_\delta^\infty P_{\theta_0}(|\hat{\theta}_{n+1} - \theta_0|^p > t) \, dt.$$

The standard derivation of the asymptotic distribution of the maximum likelihood estimator coupled with large deviation analysis yields that

$$(18) \quad \int_0^\delta P_{\theta_0}(|\hat{\theta}_{n+1} - \theta_0|^p > t) \, dt = (1 + o(1)) E|N(0,1)|^p (nI(\theta)\kappa(\theta))^{-p/2}.$$

As for the second integral in (17), let $s > 0$. It follows from (16) that there exists a constant $c_2 > 0$ independent of $\delta$ and $s$ such that

$$\int_\delta^\infty P_{\theta_0}(|\hat{\theta}_{n+1} - \theta_0|^p > t) \, dt$$

$$\leq \int_\delta^\infty P_{\theta_0}((\overline{\log X} - E_{\theta_0} \log X)^p e^{p\overline{\log X}} c_1^p > t)\, dt$$

$$= \int_\delta^\infty P_{\theta_0}((\overline{\log X} - E_{\theta_0} \log X)^p e^{p[(\overline{\log X} - E_{\theta_0} \overline{\log X}) + E_{\theta_0} \overline{\log X}]} c_1^p > t)\, dt$$

$$\leq \int_\delta^\infty P_{\theta_0}(\overline{\log X} > E_{\theta_0} \log X + c_2 \log t)\, dt$$

$$= \int_\delta^\infty P_{\theta_0}\left(\left(\prod_{i=1}^n X_i\right)^s > e^{(E_{\theta_0} \log X + c_2 \log t)ns}\right) dt$$

$$\leq \int_\delta^\infty ((E_{\theta_0} X^s)/e^{(sE_{\theta_0} \log X)})^n t^{-c_2 ns}\, dt$$

$$= ((E_{\theta_0} X^s)/e^{sE_{\theta_0} \log X})^n \frac{\delta^{-(c_2 ns - 1)}}{c_2 ns - 1}.$$

Combining this with (17) and (18), after setting $\delta$ large enough so that
$$(E_{\theta_0} X^s)/[\delta^{c_2 s} e^{sE_{\theta_0} \log X}] < 1$$
yields
$$E_{\theta_0}|\theta_n - \theta_0|^p = (1 + o(1))E|N(0,1)|^p (nI(\theta)\kappa(\theta))^{-p/2}.$$

Thus the first two conditions of Theorem 3 (with $p = 2$ and $p = 4$) are satisfied.

It is left to show that for some $c > 0$

(19) $$\sum_{n=1,\ldots,\infty} E_{\theta_0}[(\log((\hat{\theta_n})^{-1}))^+ \mathbb{1}(\hat{\theta_n} < c)] < \infty.$$

Since $\zeta(y) = -(1 + o(1))/y$ as $y \to 0$, it follows that there exists a constant $0 < c^* < 1 \wedge \theta_0$ such that, for any $0 < c < c^*$, $\overline{\log X} < 0$ on $\{\hat{\theta}_{n+1} < c\}$, and for any such $c$ there exists a constant $c_3 > 0$ such that $\hat{\theta}_{n+1} = \zeta^{-1}(\overline{\log X}) \geq c_3/(-\overline{\log X})$ on $\{\hat{\theta}_{n+1} < c\}$. The closer to zero that one chooses $c$, the closer to 1 one can set $c_3$. Choose $0 < s < \theta_0$ and $c$ and define $c_4$ to be such that $c_3/c > [1 + \Gamma(\theta_0 - s)/\Gamma(\theta_0)]^{1/s} = c_4$. Now

$$E_{\theta_0}[(\log((\hat{\theta}_{n+1})^{-1}))^+ \mathbb{1}(\hat{\theta}_{n+1} < c)]$$

(20) $$\leq E_{\theta_0}\{[(-\log c_3)^+ + (\log(\overline{-\log X}))^+] \mathbb{1}(\hat{\theta}_{n+1} < c)\}$$

$$\leq (-\log c_3)^+ P_{\theta_0}(\hat{\theta}_{n+1} < c) + E_{\theta_0}\{\log(\overline{-\log X})^+ \mathbb{1}(\hat{\theta}_{n+1} < c)\}.$$

Since $c < \theta_0$, standard considerations of large deviation analysis yield that $P_{\theta_0}(\hat{\theta}_{n+1} < c)$ is exponentially bounded, so that $\sum_{n=1,\ldots,\infty} P_{\theta_0}(\hat{\theta}_{n+1} < c) < \infty$. As for the last term in (20), recall that $0 < s < \theta_0$ and note that

$$E_{\theta_0}\{\log(\overline{-\log X})^+ \mathbb{1}(\hat{\theta}_{n+1} < c)\}$$

$$\leq \int_{\log c_4}^\infty P_{\theta_0}\{\log(\overline{-\log X}) > t\}\, dt + (\log c_4) P_{\theta_0}(\hat{\theta}_{n+1} < c).$$

The last term on the right-hand side sums to a finite result by the same argument just given, and the integral equals

$$\int_{\log c_4}^{\infty} P_{\theta_0}\{\overline{-\log X} > e^t\} dt$$

$$= \int_{\log c_4}^{\infty} P_{\theta_0}\left\{\exp\left(\sum_{i=1,\ldots,n} s\log(1/X_i)\right) > \exp(sne^t)\right\} dt$$

$$\leq \int_{\log c_4}^{\infty} \frac{(E_{\theta_0}\exp(s\log(1/X)))^n}{\exp(sne^t)} dt$$

$$\leq \left(\frac{\Gamma(\theta_0 - s)}{\Gamma(\theta_0)}\right)^n \int_{\log c_4}^{\infty} e^{-sne^t} dt$$

$$< \left(\frac{\Gamma(\theta_0 - s)}{\Gamma(\theta_0)}\right)^n \int_{\log c_4}^{\infty} e^{-snt} dt$$

$$= \left(\frac{\Gamma(\theta_0 - s)}{\Gamma(\theta_0)}\right)^n \frac{c_4^{-sn}}{sn},$$

from which (19) follows using the definition of $c_4$.

**A.3. Description of Monte Carlo.** Let $\psi(b) = E\exp\{-(\log(\Lambda_{\tau_b}) - b)\}$.

We wish to estimate $\gamma = \lim_{b\to\infty}\psi(b)$ by simulating $\exp\{-(\log(\Lambda_{\tau_b}) - b)\}$ $r$ times for a large $b$ and averaging the results. It is obviously efficient to use the same simulation runs to estimate $\psi(b)$ for a chosen set of values $b_1, \ldots, b_k$. We want to choose large $b_i$'s, large enough so that the simulation results $\psi(b_i)$ are approximately equal, in which case it is reasonable to assume that they are close to the limiting value $\gamma$. The accuracy of simulation results is improved by "averaging over intervals of $b$-values": define

$$\overline{\psi}(B_0, B_1) = \int_{B_0}^{B_1} \psi(b)\,db/(B_1 - B_0)$$

$$= \int_{B_0}^{B_1} \mathrm{avg}(\exp\{-\xi(b)\})\,db/(B_1 - B_0),$$

where $\xi(b)$ denotes the overshoot of $\{\log \Lambda_n\}$ over $b$, and "avg" denotes the (sample) mean of the results for the $r$ simulation runs.

Interchanging the operations of integration and averaging,

$$(B_1 - B_0)\overline{\psi}(B_0, B_1) = \mathrm{avg}\left(\int_{B_0}^{B_1}(\exp\{-\xi(b)\})\,db\right).$$

Consider the ladder variables (successive "record values") of the process $\{\log \Lambda_n\}$, and define for a given run $\Gamma_0 \equiv B_0 < \Gamma_1 < \Gamma_2 < \cdots < \Gamma_{m-1} \leq B_1 <$

$\Gamma_m$ as the $m \geq 1$ ladder variables in $(B_0, B_1]$ and the first one overshooting $B_1$. Then

$$\int_{B_0}^{B_1} (\exp\{-\xi(b)\})\, db = \sum_{i=1}^{m} \int_{a_{i-1}}^{a_i} \exp\{-(\Gamma_i - b)\}\, db,$$

where $a_0 = \Gamma_0 = B_0$, $a_1 = \Gamma_1, \ldots, a_{m-1} = \Gamma_{m-1}$ and $a_m = B_1$. It is easy to calculate the integrals, yielding

$$(B_1 - B_0)\overline{\psi}(B_0, B_1) = \mathrm{avg}\bigg(\sum_{i=1}^{m}(1 - \exp\{\Gamma_{i-1} - \Gamma_i\}) + \exp\{B_1 - \Gamma_m\} - 1\bigg).$$

This formula is easily applied by accumulating on each simulation run the terms coming from the "ladder steps," $\Gamma_i - \Gamma_{i-1}$, and using also the values $\Gamma_m - B_1$ when $B_1$ is first exceeded, ending the run.

For each of the cases $t = 0, 0.5, 1$, three simulations of $r = 5000$ runs each were performed, using $[B_0, B_1] = [10, 15], [15, 20]$ and $[20, 25]$. The results, shown in Table 1, indicate that $b \geq 10$ (corresponding to $\alpha < e^{-10} \approx 1/22{,}000$) is large enough in the present example to provide a stable estimate of $\gamma$.

The simulation runs were truncated after $n\max = 50{,}000$ observations when $[B_0, B_1] = [10, 15]$, and after 75,000 and 100,000 observations in the other two cases. In 1–2% of the runs, the boundary $B_1$ was not reached before truncation (due to $\theta$ estimates close to the null value, $\theta_0 = 1$). In most of these instances, $B_0$ was not reached either. In the latter cases, the results for those runs were divided not by $B_1 - B_0$ but by "the largest observed ladder variable"—$B_0$. When $B_0$ was not reached, the value 1 was used as the output of the run in computing the averages over the $r$ runs. Both of these adjustments seem appropriate and have a small effect on the tabulated results, presumably causing a very slight positive bias of the estimates of $\gamma$, much smaller than their standard errors.

The simulations reported in Table 3 were speeded up using linear interpolation in a table of 30,000 values of the maximum likelihood estimator over the range $[-10, 10]$ for the average of the $\log X$'s, a tactic which should not be needed when the SRRS procedure is applied to a single observed sequence.

**Acknowledgments.** The authors are grateful to the referees and especially to an Associate Editor who gave us one of the most thorough, accurate and helpful reports we have ever received.

Department of Mathematics
California Institute of Technology
Pasadena, California 91125
USA
e-mail: glorden@caltech.edu

Department of Statistics
Hebrew University of Jerusalem
91905 Jerusalem
Israel
e-mail: msmp@mscc.huji.ac.il